\documentclass[10pt]{amsart}
\usepackage{enumerate}
\usepackage{amsmath,amssymb}
\newcommand{\doublespace}
   {\addtolength{\baselineskip}{0.15\baselineskip}}

\newtheorem{pdef}{Definition}[section] %
\newtheorem{thm}[pdef]{Theorem}        

\newtheorem{lem}[pdef]{Lemma}
\newtheorem{prop}[pdef]{Proposition}

\newtheorem{remark}[pdef]{Remark}
\newcounter{equationnumber}
\renewcommand{\theequation}{\thesection.\arabic{equation}}
\def\mathletters{
    \addtocounter{equation}{1}
    \edef\@currentlabel{\theequation}
    \setcounter{equationnumber}{\value{equation}}
    \setcounter{equation}{0}
    \edef\theequation{\@currentlabel\noexpand\alph{equation}}
    }

\title{On the support of measures in multiplicative free convolution semigroups}
\author{Hao-Wei Huang and Ping Zhong}

\address{Department of Mathematics, Indiana University, 831 East 3rd Street,
Bloomington, IN 47405} \email{huang39@indiana.edu,
pzhong@indiana.edu}

\begin{document}
\maketitle \doublespace \pagestyle{myheadings} \thispagestyle{plain}
\markboth{H.-W. Huang and P. Zhong}{Project 2}

\begin{abstract} In this paper, we study the supports of measures in multiplicative free semigroups on
the positive real line and on the unit circle. We provide formulas
for the density of the absolutely continuous parts of measures in
these semigroups. The descriptions rely on the characterizations of
the images of the upper half-plane and the unit disc under certain
subordination functions. These subordination functions are
$\eta$-transforms of infinitely divisible measures with respect to
multiplicative free convolution. The characterizations also help us
study the regularity properties of these measures. One of the main
results is that the number of components in the support of measures
in the semigroups is a decreasing function of the semigroup
parameter.
\end{abstract}
\footnotetext[1]{{\it 2000 Mathematics Subject Classification:}\,
Primary 46L54, Secondary 30A99.} \footnotetext[2]{{\it Key words and
phrases.}\, Free convolution power, Cauchy transform, support,
regularity.}

\section{Introduction}

Denote by $\mathcal{M}_{\mathbb{R}_+}$ and $\mathcal{M}_\mathbb{T}$
the set of Borel probability measures on the positive real line
$\mathbb{R}_+=[0,\infty)$ and on the unit circle
$\mathbb{T}=\{z\in\mathbb{C}:|z|=1\}$, respectively. Further let
$\mathcal{M}_{\mathbb{R}_+}^\times=\mathcal{M}_{\mathbb{R}_+}\backslash\{\delta_0\}$,
and let $\mathcal{M}_\mathbb{T}^\times$ be the subset of
$\mathcal{M}_{\mathbb{T}}$ consisting of measures whose first
moments are nonzero and $\eta$-transforms (see section 2 for
definition) never vanish throughout $\mathbb{D}\backslash\{0\}$,
where $\mathbb{D}=\{z\in\mathbb{C}:|z|<1\}$. For measures $\mu$ and
$\nu$ both in $\mathcal{M}_{\mathbb{R}_+}^\times$ or in
$\mathcal{M}_\mathbb{T}^\times$, their multiplicative free
convolution, denoted by $\mu\boxtimes\nu$, can be characterized via
the $\Sigma$-transforms (see section 2 for definition). That is, the
identity $\Sigma_{\mu\boxtimes\nu}(z)=\Sigma_\mu(z)\Sigma_\nu(z)$
holds for $z$ in some appropriate region. We refer the reader to
[\ref{HV2},\ref{V4}] for more details.

One of the significant properties of multiplicative free convolution
is the existence of subordination functions. More precisely, the
$\eta$-transform $\eta_{\mu\boxtimes\nu}$ of measures $\mu$ and
$\nu$ either both in $\mathcal{M}_{\mathbb{R}_+}$ or in
$\mathcal{M}_\mathbb{T}$ is subordinated to $\eta_\mu$ in the sense
that $\eta_{\mu\boxtimes\nu}=\eta_\mu\circ\omega$ for some unique
analytic function $\omega$. The function $\omega$ is a self-mapping
of $\mathbb{C}\backslash\mathbb{R}_+$ if
$\mu,\nu\in\mathcal{M}_{\mathbb{R}_+}$ while it is a self-mapping of
$\mathbb{D}$ if $\mu,\nu\in\mathcal{M}_\mathbb{T}$
[\ref{BB3},\ref{Biane3},\ref{V2}]. One of the applications of
subordination functions is the study the regularity of
multiplicative free convolution. This fact, first noted by
Voiculescu [\ref{V2}], has been widely used in free probability
theory. For instance, if $\lambda_t$, $t\geq0$, is the
multiplicative free Brownian motion on $\mathbb{R}_+$ or
$\mathbb{T}$ then
$\eta_{\mu\boxtimes\lambda_t}=\eta_\mu\circ\omega_t$ for some
analytic function $\omega_t$ whose $\eta$-transform is a
$\boxtimes$-infinitely divisible measure. It turns out that the
measure $\mu\boxtimes\lambda_t$ is absolutely continuous for any
$t>0$, and its density can be described in terms of $\omega_t$. More
importantly, the number of components in the support of
$\mu\boxtimes\lambda_t$ is a decreasing function of $t$. We refer
the reader to [\ref{Biane2},\ref{P2}] for more details. The same
conclusion also holds for additive free convolution [\ref{Biane1}].

For $n\in\mathbb{N}$, and $\mu$ in $\mathcal{M}_{\mathbb{R}_+}$ or
$\mathcal{M}_\mathbb{T}$, the $n$-fold multiplicative free
convolution $\mu\boxtimes\cdots\boxtimes\mu$ is denoted by
$\mu^{\boxtimes n}$. The measure $\mu^{\boxtimes n}$ can be
characterized by its $\Sigma$-transform, that is, the identity
$\Sigma_{\mu^{\boxtimes n}}=\Sigma_\mu^n$ holds in some appropriate
region. This can be generalized to any convolution power $t\geq1$.
More precisely, if $\mu\in\mathcal{M}_{\mathbb{R}_+}^\times$ (resp.
$\mu\in\mathcal{M}_{\mathbb{T}}^\times$) and $t\geq1$ there there
exists a semigroup $\{\mu_t:t\geq1\}$ contained in
$\mathcal{M}_{\mathbb{R}_+}^\times$ (resp. contained in
$\mathcal{M}_{\mathbb{T}}^\times$) such that the identity
$\Sigma_{\mu_t}=\Sigma_\mu^t$ holds in some appropriate domain,
where the $t$th power is taken appropriately. The measure $\mu_t$
coincides with $\mu^{\boxtimes n}$ if $t=n$ is an integer. With the
help of the subordination in this context, it was shown that the
measure $\mu_t$ has no singular continuous part and its density is
analytic wherever it is positive [\ref{BB1},\ref{BB2}].

In this paper, we study the supports of measures in the semigroup
$\{\mu_t:t\geq1\}$ associated with the measure $\mu$ which in
$\mathcal{M}_{\mathbb{R}_+}^\times$ or in
$\mathcal{M}_\mathbb{T}^\times$. The main tools used in this study
are the properties of subordination functions established in
[\ref{BB2}]. By the methods developed in this paper, we show that
the number of components in the support of $\mu_t$ is a decreasing
function of $t$. The corresponding subordination functions are shown
to be the $\eta$-transforms of $\boxtimes$-infinitely divisible
measures and their ranges are also analyzed. Another purpose of this
paper is to give an implicit formula for the density of the
absolutely continuous part of $\mu_t$.

The paper is organized as follows. Section 2 contains some
definitions and preliminaries. Section 3 and Section 4 investigate
the topological properties of measures in the semigroups associated
with measures $\mu$ on the positive real line and the unit circle,
respectively.

\section{Preliminary}

Following the notation in the introduction, the $\psi$-transform of
$\mu$ is defined as
\[\psi_\mu(z)=\int\frac{zs}{1-zs}\;d\mu(s),\] which
is analytic on $\Omega=\mathbb{C}\backslash[0,+\infty)$ if
$\mu\in\mathcal{M}_{\mathbb{R}_+}$ and analytic on $\mathbb{D}$ if
$\mu\in\mathcal{M}_\mathbb{T}$. The $\eta$-transform of $\mu$ is
defined as
\[\eta_\mu=\frac{\psi_\mu}{1+\psi_\mu}\] on the same domain as the
$\psi$-transform. The analytic way to obtain the multiplicative free
convolution is by using the inverse of the $\eta$-transform. A
measure $\mu$ is said to be $\boxtimes$-infinitely divisible if for
any integer $n$ there exists a measure $\mu_n$ such that
\[\mu=\underbrace{\mu_n\boxtimes\cdots\boxtimes\mu_n}_{n\;\;\mathrm{factors}}.\] We
refer the reader to [\ref{HV1},\ref{HV2}] for more information about
multiplicative free convolution and $\boxtimes$-infinite
divisibility. In the following two subsections, we briefly introduce
some properties of these transforms and the multiplicative free
semigroups associated with a measure in
$\mathcal{M}_{\mathbb{R}_+}^\times$ or
$\mathcal{M}_\mathbb{T}^\times$.

\subsection{Multiplicative Free Convolution Semigroup on $\mathbb{R}_+$}

Measures in $\mathcal{M}_{\mathbb{R}_+}^\times$ can be characterized
by their $\eta$-transforms, which is stated in the following
proposition.

\begin{prop} \label{chareta1} Let
$\eta:\Omega\to\mathbb{C}\setminus\{0\}$ be an analytic function
such that $\eta(\overline{z})=\overline{\eta(z)}$ for all
$z\in\Omega$. Then the following statements are equivalent.
\begin{enumerate} [$\qquad(1)$]
\item {There exists a measure $\mu\in\mathcal{M}_{\mathbb{R}_+}^\times$ such that
$\eta=\eta_\mu$.}
\item {The function $\eta$ satisfies $\eta(0-)=0$ and $\arg\eta(z)\in[\arg z,\pi)$ for $z\in\mathbb{C}^+$.}
\end{enumerate}
\end{prop}

Any measure $\mu\in\mathcal{M}_{\mathbb{R}_+}$ can be recovered from
its $\eta$-transform by Stieltjes inversion formula. Indeed, the
identity
\[G_\mu\left(\frac{1}{z}\right)=\frac{z}{1-\eta_\mu(z)},\;\;\;\;\;z\in\Omega,\]
where $G_\mu$ is the Cauchy transform of $\mu$, shows that the
family of measures $\{\mu_\epsilon\}_{\epsilon>0}$ defined as
\[d\mu_\epsilon(1/x)=\frac{1}{\pi}\Im\left(\frac{x+i\epsilon}{1-\eta_\mu(x+i\epsilon)}\right)dx,\;\;\;\;\;x>0,\]
converges to $\nu$ weakly as $\epsilon\to0$, where
$d\nu(x)=d\mu(1/x)$. Note that $\nu$ is not always in
$\mathcal{M}_{\mathbb{R}_+}$ since $\nu(\mathbb{R}_+)=1-\mu(\{0\})$.

If $\mu\in\mathcal{M}_{\mathbb{R}_+}^\times$ then $\eta_\mu'(z)>0$
for $z<0$, and therefore $\eta_\mu|(-\infty,0)$ is invertible. Let
$\eta_\mu^{-1}$ be the inverse of $\eta_\mu$ and set
\[\Sigma_\mu(z)=\frac{\eta_\mu^{-1}(z)}{z},\] where $z<0$ is sufficiently small. For
$n\in\mathbb{N}$, the multiplicative free convolution power
$\mu^{\boxtimes n}$ of $\mu$ satisfies the identity
\[\Sigma_{\mu^{\boxtimes n}}(z)=\Sigma_\mu^n(z),\] where $z<0$ is sufficiently small.
The generalization to any multiplicative free convolution power
greater than one is stated below [\ref{BB2}].

\begin{thm} \label{Phit1prop}
If $\mu\in\mathcal{M}_{\mathbb{R}_+}^\times$ and $t>1$ then the
following statements hold.
\begin{enumerate} [$\qquad(1)$]
\item {There exists a unique measure $\mu^{\boxtimes t}\in\mathcal{M}_{\mathbb{R}_+}^\times$ such that
$\Sigma_{\mu^{\boxtimes t}}(z)=\Sigma_\mu^t(z)$ for $z<0$
sufficiently close to zero.}
\item {There exists an analytic function $\omega_t:\Omega\to\Omega$ such that $\omega_t((-\infty,0))\subset(-\infty,0)$,
$\omega_t(0-)=0$, $\arg\omega_t(z)\in[\arg z,\pi)$ for
$z\in\mathbb{C}^+$, and $\eta_{\mu^{\boxtimes
t}}=\eta_\mu\circ\omega_t$.}
\item {The function $\Phi_t:\Omega\to\mathbb{C}\backslash\{0\}$
defined by
\[\Phi_t(z)=z\left[\frac{z}{\eta_\mu(z)}\right]^{t-1},\;\;\;\;\;z\in\Omega,\]
where the power is taken to positive for $z<0$, satisfies
$\Phi_t(\omega_t(z))=z$ for $z\in\Omega$.}
\item {Let $x>0$. Then the point $1/x$ is an atom of $\mu^{\boxtimes t}$ if
$\mu(\{x^{-1/t}\})>(t-1)/t$, in which case we have
\[\mu^{\boxtimes t}(1/x)=t\mu(\{x^{-1/t}\})-(t-1).\]}
\end{enumerate}
\end{thm}

Next, we introduce some special mappings and sets which connect
Theorem \ref{Phit1prop} and the global inversion theorem.

Denote by $\mathcal{S}_\pi$ the strip $\{z\in\mathbb{C}:|\Im
z|<\pi\}$. Further, let
$\mathcal{S}_\pi^+=\mathcal{S}_\pi\cap\mathbb{C}^+$ and
$\mathcal{S}_\pi^-=\mathcal{S}_\pi\cap\mathbb{C}^-$. The map
$\Lambda(z)=-e^z:\mathbb{C}\to\mathbb{C}\backslash\{0\}$ is a
conformal mapping from $\mathcal{S}_\pi$ onto $\Omega$.
Particularly, we have $\Lambda(\mathcal{S}_\pi^+)=\mathbb{C}^-$ and
$\Lambda(\mathcal{S}_\pi^-)=\mathbb{C}^+$. If $\Lambda^{-1}$ is the
inverse of this conformal mapping then
\begin{equation} \label{argLambda1}
\Im\Lambda^{-1}(z)=-\pi+\arg z,\;\;\;\;\;z\in\mathbb{C}^+
\end{equation} and
\begin{equation} \label{argLambda2}
\Im\Lambda^{-1}(z)=\pi+\arg z,\;\;\;\;\;z\in\mathbb{C}^-.
\end{equation}
For any measure $\mu\in\mathcal{M}_{\mathbb{R}_+}^\times$, by
Proposition \ref{chareta1} we have
\begin{equation} \label{argeta1}
\arg[\eta_\mu(\Lambda(z))]\in(-\pi,-\pi+\Im
z],\;\;\;\;\;z\in\mathcal{S}_\pi^+
\end{equation} and
\begin{equation} \label{argeta2}
\arg[\eta_\mu(\Lambda(z))]\in[\pi+\Im
z,\pi),\;\;\;\;\;z\in\mathcal{S}_\pi^-.
\end{equation} Moreover, the function
\[\kappa_\mu(z)=\frac{z}{\eta_\mu(z)}\] is analytic on
$\Omega$ since $\eta_\mu$ never vanishes on $\Omega$. For any $t>1$,
let
\begin{equation} \label{Htdef}
H_t(z)=z+(t-1)l_\mu(z),\;\;\;\;\;z\in\mathcal{S}_\pi,
\end{equation}
where
\begin{equation} \label{lt}
l_\mu(z)=\Lambda^{-1}[-\kappa_\mu(\Lambda(z))].
\end{equation}
Next, observe that
\begin{equation} \label{upperH}
\Im H_t(z)\in[\Im z,t\Im z),\;\;\;\;z\in\mathcal{S}_\pi^+
\end{equation} or, equivalently, $\Im l_t(z)\in[0,\Im z)$, $z\in\mathcal{S}_\pi^+$.
Indeed, if $z\in\mathcal{S}_\pi^+$ then
\[\arg\kappa_\mu(\Lambda(z))=\arg\Lambda(z)-\arg\eta_\mu(\Lambda(z))\in[0,\Im
z)\] by (\ref{argeta1}), which yields
$\arg[-\kappa_\mu(\Lambda(z))]\in[-\pi,-\pi+\Im z)$ and
$-\kappa_\mu(\Lambda(z))\in\mathbb{C}^-$. Hence by
(\ref{argLambda2}) we have
\[\Im\Lambda^{-1}[-\kappa_\mu(\Lambda(z))]=\pi+\arg[-\kappa_\mu(\Lambda(z))]\in[0,\Im
z),\] as desired.

The following theorem, obtained by choosing $h=\pi$ and $k=t\pi$ in
[Theorem 4.9, \ref{BB2}], plays an important role in the
investigation of the support of $\mu^{\boxtimes t}$.

\begin{thm} \label{Htprop}
If $H_t$ is the analytic function defined in $(\ref{Htdef})$ then
the following statements hold.
\begin{enumerate} [$\qquad(1)$]
\item {There exists an analytic function $\varpi_t:\mathcal{S}_\pi\to\mathcal{S}_\pi$ such that
$\varpi_t$ extends continuously to $\overline{\mathcal{S}}_\pi$,
$|\Im\varpi_t(z)|\leq|\Im z|$,
$\varpi_t(\overline{z})=\overline{\varpi(z)}$, and
$H_t(\varpi_t(z))=z$ for all $z\in\mathcal{S}_\pi$.}
\item {The function $\varpi_t$ satisfies
\[\frac{|z_1-z_2|}{2(t+1)}\leq|\varpi_t(z_1)-\varpi_t(z_2)|,\;\;\;\;\;z_1,z_2\in\overline{\mathcal{S}}_\pi.\]}
\item {The set $\{z\in\mathcal{S}_\pi:H_t(z)\in\mathcal{S}_\pi\}$ coincides with $\varpi_t(\mathcal{S}_\pi)$
and is a simply connected domain whose boundary consists of two
simple curves, $\varpi_t(\mathbb{R}\pm i\pi)$.}
\item {If $\alpha\in\partial\mathcal{S}_\pi$ and $\varpi_t(\alpha)\in\mathcal{S}_\pi$ then
$\varpi_t$ can be continued analytically to a neighborhood of
$\alpha$.}
\end{enumerate}
\end{thm}

\subsection{Multiplicative Free Convolution Semigroup on $\mathbb{T}$}

The following proposition characterizes functions which are
$\eta$-transforms.

\begin{prop} \label{chareta2}
If $\eta:\mathbb{D}\to\mathbb{C}$ is analytic then the following conditions are equivalent.
\begin{enumerate} [\qquad$(1)$]
\item {There exists a measure $\mu\in\mathcal{M}_\mathbb{T}$ such that $\eta=\eta_\mu$.}
\item {We have $\eta(0)=0$ and $|\eta(z)|<1$ for all $z\in\mathbb{D}$.}
\item {We have $|\eta(z)|\leq|z|$ for all $z\in\mathbb{D}$.}
\end{enumerate}
\end{prop}

Any measure $\mu\in\mathcal{M}_\mathbb{T}$ can be recovered from its
$\eta$-transform. Indeed, the identity
\[\frac{1}{2\pi}\frac{1+\eta_\mu(z)}{1-\eta_\mu(z)}=\frac{1}{2\pi}\int_{-\pi}^\pi
\frac{\zeta+z}{\zeta-z}\;d\mu(1/\zeta),\] whose real part is the
Poisson integral of the measure $d\mu(1/\zeta)$ indicates that the
family of measures $\{\mu_\varepsilon\}_{\varepsilon>0}$ defined as
\[d\mu_\varepsilon(e^{i\theta})=\frac{1}{2\pi}\frac{1-|\eta_\mu(\varepsilon e^{i\theta})|^2}
{|1-\eta_\mu(\varepsilon e^{i\theta})|^2}\;d\theta\] converges to
$\nu$ weakly on $\mathbb{T}$ as $\varepsilon\downarrow0$, where
$d\mu(\zeta)=d\mu(1/\zeta)$.

Recall that measures in $\mathcal{M}_\mathbb{T}^\times$ have nonzero
mean. That is, if $\mu\in\mathcal{M}_\mathbb{T}^\times$ then
\[\eta_\mu'(0)=\int_\mathbb{T}\zeta\;d\mu(\zeta)\neq0,\] and
therefore $\eta_\mu$ is invertible in a neighborhood of zero. This
shows that the inverse $\eta_\mu^{-1}$ is defined for sufficiently
small $z$, and so is
\[\Sigma_\mu(z)=\frac{\eta_\mu^{-1}(z)}{z}.\] For
$n\in\mathbb{N}$, the multiplicative free convolution power
$\mu^{\boxtimes n}$ satisfies
\[\Sigma_{\mu^{\boxtimes n}}(z)=\Sigma_\mu^n(z)\] in a neighborhood of zero.
The following theorem is the generalization to and multiplicative
free convolution power $t\geq1$ [\ref{BB2}].

\begin{thm} \label{Phit2prop}
If $\mu\in\mathcal{M}_\mathbb{T}^\times$ and $t>1$ then the
following statements hold.
\begin{enumerate} [$\qquad(1)$]
\item {There exists a measure $\mu_t\in\mathcal{M}_\mathbb{T}^\times$
such that $\Sigma_{\mu_t}(z)=\Sigma_\mu^t(z)$ holds in a
neighborhood of zero.}
\item {There exists an analytic function
$\omega_t:\mathbb{D}\to\mathbb{D}$ such that $|\omega_t(z)|\leq|z|$
and $\eta_{\mu_t}(z)=\eta_\mu(\omega_t(z))$ for all $z\in\Omega$.}
\item {The function $\omega_t$ is given by
\[\omega_t(z)=\eta_{\mu_t}(z)\left[\frac{z}{\eta_{\mu_t}(z)}\right]^{1/t},\;\;\;\;z\in\mathbb{D}.\]}
\item {The analytic function
$\Phi_t:\mathbb{D}\to\mathbb{C}$ defined by
\[\Phi_t(z)=z\left[\frac{z}{\eta_\mu(z)}\right]^{t-1},\;\;\;z\in\mathbb{D},\]
satisfies $\Phi_t(\omega_t(z))=z$ for $z\in\mathbb{D}$.}
\item {A point $1/\zeta$ is an atom of $\mu_t$ if $\mu(\{1/\omega_t(\zeta)\})>(t-1)/t$, in which
case we have
\[\mu_t(1/\zeta)=t\mu(\{1/\omega_t(\zeta)\})-(t-1).\]}
\end{enumerate}
\end{thm}

\begin{remark} \emph{For $\mu\in\mathcal{M}_\mathbb{T}^\times$, let $\kappa_\mu(z)=z/\eta_\mu(z)$,
$z\in\mathbb{D}$.
Observe that the function $\Phi_t(z)=z\kappa_\mu(z)^{t-1}$ depends
on the choice of extracting roots, and therefore the measure $\mu_t$
in Theorem \ref{Phit2prop}(1) is not unique. However, there is only
one measure $\mu_t$ satisfying $\eta_{\mu_t}=\eta_\mu\circ\omega_t$
if $\Phi_t$ is chosen.}
\end{remark}

The function $\omega_t$ in the preceding result is obtained as a
consequence of the following global inversion theorem [\ref{BB2}].

\begin{thm} \label{GIT2}
Let $\Phi:\mathbb{D}\to\mathbb{C}\cup\{\infty\}$ be a
meromorphic function such that $\Phi(0)=0$ and $|z|\leq|\Phi(z)|$
for all $z\in\mathbb{D}$. Then the following statements hold.
\begin{enumerate} [$\qquad(1)$]
\item {There exists a continuous function $\omega:\overline{\mathbb{D}}\to\overline{\mathbb{D}}$ such that
$\omega(\mathbb{D})\subset\omega(\mathbb{D})$, $\omega|\mathbb{D}$
is analytic, and $\Phi(\omega(z))=z$ for all $z\in\mathbb{D}$.}
\item {If $\zeta\in\mathbb{T}$ is such that $|\omega(\zeta)|<1$ then $\omega$ can be continued analytically
to a neighborhood of $\zeta$.}
\item {The set $\{z\in\mathbb{D}:|\Phi(z)|<1\}$ is simply connected bounded by a simple closed curve.
This set coincides with $\omega(\mathbb{D})$ and its boundary is
$\omega(\mathbb{T})$.}
\item {If
$z\in\omega(\overline{\mathbb{D}})\cap\mathbb{T}$ then the entire
radius $\{rz:0\leq r<1\}$ is contained in $\omega(\mathbb{D})$.}
\end{enumerate}
\end{thm}

\setcounter{equation}{0}
\section{Support of the measure $\mu^{\boxtimes t}$ on $\mathbb{R}_+$}

Throughout this section, we fix a measure
$\mu\in\mathcal{M}_{\mathbb{R}_+}^\times$ and investigate the
support of $\mu^{\boxtimes t}$, $t\geq1$, which is the unique
measure defined in Theorem \ref{Phit1prop}. Let
\[\Gamma_t=\{z\in\mathcal{S}_\pi^-:H_t(z)\in\mathcal{S}_\pi\},\] where $H_t$ is the
function defined as in (\ref{Htdef}). In the following proposition,
we describe the set $\Gamma_t$ in terms of $\kappa_\mu$.

\begin{prop} \label{Gammat}
The set $\Gamma_t$ can be expressed as
\begin{equation} \label{Gammat1}
\Gamma_t=\left\{z\in\mathcal{S}_\pi^-:\frac{-\arg\kappa_\mu(\Lambda(z))}{\arg\Lambda(z)}<\frac{1}{t-1}\right\}.
\end{equation} Moreover, $\Gamma_t$ is a simply connected domain whose boundary consists of
two simple curves one of which is the real line.
\end{prop}

\begin{pf} Since $H_t(\overline{z})=\overline{H_t(z)}$ for any $z\in\mathcal{S}_\pi$,
by (\ref{upperH}) we see that $\Im H_t(z)\in(t\Im z,\Im z]$ for
$z\in\mathcal{S}_\pi^-$. This shows that a point
$z\in\mathcal{S}_\pi^-$ satisfies $H_t(z)\in\mathcal{S}_\pi$ if and
only if
\[\Im l_\mu(z)>-\frac{\pi+\Im z}{t-1}=-\frac{\arg\Lambda(z)}{t-1},\] Since
$-\kappa_\mu(\Lambda(z))\in\mathbb{C}^+$ and
$\arg[-\kappa_\mu(\Lambda(z))]=\pi+\arg\kappa_\mu(\Lambda(z))$ for
$z\in\mathcal{S}_\pi^-$, the description (\ref{Gammat1}) for
$\Gamma_t$ follows from (\ref{argLambda1}). Finally, the fact
$H_t(\overline{z})=\overline{H_t(z)}$ shows that the set $\{z\in
S_\pi:H_t(z)\in S_\pi\}$ is symmetry with respect to the real line,
and hence the second assertion follows.
\end{pf} \qed

\begin{prop} \label{extension1}
The function $l_\mu$ defined in $(\ref{lt})$ has a continuous
extension to the boundary $\overline{\Gamma}_t$ and the extension is
Lipschitz continuous. Consequently, the function
$\kappa_\mu\circ\Lambda$ has a continuous extension to
$\overline{\Gamma}_t$.
\end{prop}

\begin{pf} By Theorem \ref{Htprop}(2), we have
\[\left|\frac{H_t(z_1)-H_t(z_2)}{z_1-z_2}\right|=\left|1+(t-1)\frac{l_\mu(z_1)-l_\mu(z_2)}
{z_1-z_2}\right|\leq2(t+1),\;\;\;\;\;z_1,z_2\in\Gamma_t,\] which
yields that
\[\left|\frac{l_\mu(z)-l_\mu(z)}
{z_1-z_2}\right|\leq\frac{2t+3}{t-1},\;\;\;\;\;z_1,z_2\in\Gamma_t.\]
This shows that $l_\mu$ extends continuously to
$\overline{\Gamma}_t$ whose extension is Lipschitz continuous. Since
$\kappa_\mu(\Lambda(z))=\exp[l_\mu(z)]$ for all
$z\in\mathcal{S}_\pi$, it follows that $\kappa\circ\Lambda$ extends
continuously to $\overline{\Gamma}_t$.
\end{pf} \qed

The characterization and non-vanishing of the $\eta$-transform of a
measure $\mu$ in $\mathcal{M}_{\mathbb{R}_+}^\times$ gives that
\begin{equation} \label{lower}
\arg\kappa_\mu(z)\in[-\pi+\arg z,0),\;\;\;\;\;z\in\mathbb{C}^+.
\end{equation} This implies that
\[\kappa_\mu(z)=\exp[u(z)],\;\;\;\;\;z\in\mathbb{C}^+,\]
where $u$ is an analytic function on $\Omega$ satisfying
$u(\overline{z})=\overline{u(z)}$ for $z\in\mathbb{C}^+$ and
$u(\mathbb{C}^+)\subset\mathbb{C}^-\cup\mathbb{R}$. As a consequence
the Nevanlinna representation, $u$ can be written as
\begin{equation} \label{repre1}
u(z)=a-bz+\int_0^\infty\frac{1+zs}{z-s}\;d\rho(s),\;\;\;\;\;z\in\Omega,
\end{equation}
where $a\in\mathbb{R}$, $b\geq0$, and $\rho$ is some positive Borel
measure on $[0,\infty)$. In the following lemmas we provide some
properties of the function $u$.

\begin{lem} \label{repre1b}
If $u$ has the expression $(\ref{repre1})$ then $b=0$, and therefore
\[\kappa_\mu(z)=\exp\left(a+\int_0^\infty\frac{1+zs}{z-s}\;d\rho(s)\right),\;\;\;\;\;z\in\Omega.\]
\end{lem}

\begin{pf} First observe
that
\[\lim_{x\to+\infty}(1+\psi_\mu(-x))=\lim_{x\to+\infty}\int_0^\infty\frac{d\mu(s)}{1+xs}=0\]
by dominated convergence theorem. Moreover, for $x\geq1$ we have
\[0<\int_0^\infty\frac{d\mu(s)}{x(s+1)}\leq1+\psi_\mu(-x)\]
or, equivalently,
\[0<\frac{1}{1+\psi_\mu(-x)}\leq\frac{x}{c},\] where
\[c=\int_0^\infty\frac{d\mu(s)}{s+1}\] is a finite positive number.
This implies that
\[0\leq\lim_{x\to+\infty}\frac{\log[1+\psi_\mu(-x)]}{-x}\leq\lim_{x\to+\infty}\frac{\log x-\log
c}{x}=0.\] Then the above discussions and the expression
\[\log\kappa_\mu(-x)=\log x+\log[1+\psi_\mu(-x)]-\log[-\psi_\mu(-x)],\;\;\;\;\;x>0,\]
yield that
\[\lim_{x\to-\infty}\frac{u(x)}{x}=\lim_{x\to+\infty}\frac{\log\kappa_\mu(-x)}{-x}=0.\]
Since
\begin{align*}
\lim_{x\to-\infty}\frac{\int_0^\infty\frac{1+xs}{x-s}\;d\rho(s)}{x}=
\lim_{x\to+\infty}\int_0^\infty\frac{1-xs}{x(x+s)}\;d\rho(s)=0
\end{align*} by dominated convergence theorem, we must have $b=0$, as
desired.
\end{pf} \qed

In the sequel, the measure $\rho$ will be the unique measure in the
Nevanlinna representation \ref{repre1} of $u$.

\begin{lem} \label{decreasing1}
For any fixed $r>0$, the function
\[g(r,\theta)=-\frac{\Im u(re^{i\theta})}{\theta}\]
is decreasing on $(0,\pi)$ and $\lim_{\theta\to\pi^-}g(r,\theta)=0$.
Consequently, $\arg\kappa_\mu(z)=\Im u(z)\in(-\pi+\arg z,0]$ for
$z\in\mathbb{C}^+$.
\end{lem}

\begin{pf} First observe that for $\theta\in(0,\pi)$ the function $g(r,\theta)$ can be expressed as
\[g(r,\theta)=\frac{r\sin\theta}{\theta}\int_0^\infty\frac{s^2+1}{r^2-2rs\cos\theta+s^2}\;d\rho(s)\]
by Lemma \ref{repre1b}. To show that $g(r,\cdot)$ is decreasing on
$(0,\pi)$, it suffices to show that it has a negative derivative.
This follows from the facts
\[\frac{d}{d\theta}\left(\frac{\sin\theta}{\theta}\right)=\frac{\cos\theta}{\theta^2}(\theta-\tan\theta)<0\]
and
\[\frac{d}{d\theta}\left(\frac{1}{r^2-2rs\cos\theta+s^2}\right)=\frac{-2rs\sin\theta}{(r^2-2rs\cos\theta+s^2)^2}<0\]
for any $s\in[0,\infty)$ and $\theta\in(0,\pi)$. Since
$u(x)\in\mathbb{R}$ for $x<0$, it follows that
$\lim_{\theta\to\pi^-}g(r,\theta)=0$ for any $r>0$. The last
assertion follows from the above discussion, (\ref{lower}), and the
continuity of $u$ on $\mathbb{C}^+$.
\end{pf} \qed

For $t>0$, the function $\Phi_t$ defined in Theorem
\ref{Phit1prop}(3) can be expressed as
\begin{equation} \label{Phit1def}
\Phi_t(z)=z\exp[(t-1)u(z)].
\end{equation}
Then there exists an analytic function $\omega_t:\Omega\to\Omega$
satisfying the properties listed in Theorem \ref{Phit1prop}(2).
Indeed, we have the relations
\[\Phi_t=\Lambda\circ H_t\circ\Lambda^{-1}\;\;\;\;\;\mathrm{and}\;\;\;\;\;
\omega_t=\Lambda\circ\varpi_t\circ\Lambda^{-1}.\] The subordination
function $\omega_t$ has an important property, which is stated in
the following result.

\begin{prop} The function $\omega_t$ is the $\eta$-transform of some
$\boxtimes$-infinitely divisible measure in
$\mathcal{M}_{\mathbb{R}_+}^\times$.
\end{prop}

\begin{pf} Since $\omega_t$ satisfies the conditions in Theorem
\ref{chareta1}, we must have $\omega_t=\eta_{\nu_t}$ for some
measure $\nu_t\in\mathcal{M}_{\mathbb{R}_+}^\times$. It is clear
that the $\Sigma$-transform of $\nu_t$ is given by
$\Sigma_{\nu_t}(z)=\Phi_t(z)/z=\exp[(t-1)u(z)]$. Then the desired
result follows from [Theorem 6.12, \ref{HV2}].
\end{pf} \qed

Let $\mu^{\boxtimes t}$ be the unique measure in
$\mathcal{M}_{\mathbb{R}_+}^\times$ such that
\[\eta_{\mu^{\boxtimes}}=\eta_\mu\circ\omega_t.\] Our analysis of the support of $\mu^{\boxtimes t}$ will be based on
the functions $g:(0,\infty)\to\mathbb{R}^+\cup\{+\infty\}$ and
$A_t:\mathbb{R}^+\to[0,\pi]$ which are defined as
\[g(r)=\int_0^\infty\frac{r(s^2+1)}{(r-s)^2}\;d\rho(s)\] and
\[A_t(r)=\inf\left\{\theta\in(0,\pi):\frac{-\Im
u(re^{i\theta})}{\theta}<\frac{1}{t-1}\right\},\] respectively. The
following set, associated with the function $g$,
\[V_t^+=\left\{r\in(0,\infty):g(r)>\frac{1}{t-1}\right\}\] will also
play an important role in the investigation of the support of
$\mu^{\boxtimes t}$. Let
\[\Omega_t=\Lambda(\Gamma_t).\]
The following lemma provides some basic properties about the
functions and set defined above.

\begin{lem} \label{property1}
Let $t>1$. Then we have
\begin{equation} \label{Omegat1}
\Omega_t=\{re^{i\theta}:A_t(r)<\theta<\pi\;\;\mathrm{and}\;\;r\in(0,\infty)\},
\end{equation}
\begin{equation} \label{partialomegat1}
\partial\Omega_t=\{re^{iA_t(r)}:r\in(0,\infty)\}\cup(-\infty,0].
\end{equation} For any $r>0$, we have
\begin{equation} \label{Atr}
A_t(r)\in[0,\pi)
\end{equation} and
\begin{equation} \label{Vt+1}
V_t^+=\{r\in(0,\infty):A_t(r)>0\}.
\end{equation} Moreover, for any
$r\in(0,\infty)$ we have
\begin{equation} \label{lim}
\lim_{\theta\downarrow A_t(r)}\frac{-\Im
u(re^{i\theta)}}{\theta}\leq\frac{1}{t-1},
\end{equation} where the equality holds
if $r\in V_t^+$, that is,
\begin{equation} \label{Imu}
\Im u(re^{iA_t(r)})=-\frac{A_t(r)}{t-1},\;\;\;\;\;r\in V_t^+.
\end{equation}
\end{lem}

\begin{pf} By Proposition
\ref{Gammat} and Lemma \ref{decreasing1} we see that
$re^{i\theta}\in\Omega_t$ if and only if $\theta\in(0,\pi)$ and
$g(r,\theta)<1/(t-1)$. Since $g(r,\theta)$ is a decreasing function
on $(0,\pi)$ for $r>0$, it is clear that (\ref{Omegat1}) and
(\ref{partialomegat1}) hold by the definition of $A_t(t)$. Moreover,
since $\lim_{\theta\to\pi^-}g(r,\theta)=0$ for any $r>0$, $A_t(r)$
must belong to the interval $[0,\pi)$. Next, observe that for any
$r>1$,
\[\lim_{\theta\to0^+}g(r,\theta)=g(r)\] by the monotone convergence
theorem, and therefore the equation (\ref{Vt+1}) holds. The
inequality in (\ref{lim}) follows from the above discussion. If the
strict inequality in (\ref{lim}) occurs for some $r\in V_t^+$ then
it will violate the definition of $A_t(r)$, whence the proof is
complete.
\end{pf} \qed

\begin{prop} For any $re^{i\theta}\in\Omega_t$, the arc
$\{re^{i\phi}:\theta\leq\phi<\pi\}$ is contained in $\Omega_t$.
Consequently, the set $\Omega_t$ consists of one connected component
and $\partial\Omega_t$ is a simple curve.
\end{prop}

\begin{prop} \label{extension1}
The function $u(z)$ has a continuous extension to
$\overline{\Omega}_t$. Consequently, $\Phi_t$ and $\omega_t$ extend
continuously to $\overline{\Omega}_t$ and
$\mathbb{C}^+\cup\mathbb{R}$, respectively.
\end{prop}

\begin{pf} Since $u\circ\Lambda$ has a continuous extension to
$\overline{\Gamma}_t$ by Proposition \ref{extension1} and
$\overline{\Omega}_t=\Lambda(\overline{\Gamma}_t)$, the desired
result follows.
\end{pf} \qed

\begin{lem} \label{sconvex1}
If $g$ is bounded on some open interval $I$ then $\rho(I)=0$ and $g$
is strictly convex on $I$. In particular, this is true if $I$ is
contained in $(V_t^+)^c$.
\end{lem}

\begin{pf} Suppose that $g$ is bounded by
$M$ on $I$. By countable additivity of $\rho$, it suffices to show
that $\rho(J)=0$ for any closed interval $J$ contained in $I$. Let
$c=\min\{x:x\in J\}$ and $[a,b]\subset J$. Then $g((a+b)/2)\leq M$
gives
\[M\geq\int_a^b\frac{r(s^2+1)}{(r-s)^2}\;d\rho(s)\geq\int_a^b\frac{c(c^2+1)}{\left(\frac{b-a}{2}\right)^2}\;d\rho(s)
=4c(c^2+1)\frac{\rho((a,b))}{(b-a)^2}\] or, equivalently,
\[\frac{\rho((a,b))}{(b-a)^2}\leq\frac{M}{4c(c^2+1)}<\infty\] since $c\neq0$. Since
$[a,b]$ can be an arbitrary subinterval in $J$, we conclude that
$\rho(J)=0$, as desired.
\end{pf} \qed

Observe that the mapping $r\mapsto re^{iA_t(r)}$ is a homeomorphism
of $(0,\infty)$ onto
$\partial\Omega_t\cap(\mathbb{C}^+\cup(0,+\infty))$. It turns out
that the function $h_t:(0,\infty)\to(0,\infty)$ defined as
\[h_t(r)=\Phi_t(re^{iA_t(r)})\]
is a homeomorphism of $(0,\infty)$. We are now at the position to
introduce the main theorem of this section. For any measure $\nu$,
denote by $\nu^{\mathrm{ac}}$ and $\mathrm{supp}(\nu)$ the
absolutely continuous part and support of $\nu$, respectively.

\begin{thm} \label{density1}
Suppose that $\mu$ is a measure in
$\mathcal{M}_{\mathbb{R}_+}^\times$ and $t>1$. Let
\[S_t=\left\{\frac{1}{h_t(r)}:r\in V_t^+\right\}.\] Then the following statements hold.
\begin{enumerate} [$\qquad(1)$]
\item {The measure $(\mu^{\boxtimes t})^{\mathrm{ac}}$ is concentrated on the closure of $S_t$.}
\item {The density of $(\mu^{\boxtimes t})^{\mathrm{ac}}$ is analytic on the set
$S_t$ and is given by
\[\frac{d(\mu^{\boxtimes t})^{\mathrm{ac}}}{dx}\left(\frac{1}{h_t(r)}\right)=
\frac{1}{\pi}\frac{h_t(r)l_t(r)\sin\theta_t(r)}{1-2l_t(r)\cos\theta_t(r)+l_t^2(r)},\;\;\;\;\;r\in
V_t^+,\] where
\[l_t(r)=r\exp\Re u(re^{iA_t(r)})\] and
\[\theta_t(r)=\frac{tA_t(r)}{t-1}\] for $r\in V_t^+$.}
\item {The number of components in $\mathrm{supp}(\mu^{\boxtimes t})^{\mathrm{ac}}$ is a decreasing function of $t$.}
\end{enumerate}
\end{thm}

\begin{pf} Let $z=re^{iA_t(r)}$, $r>0$. First claim that
$\Im\eta_\mu(z)=0$ if and only if $r\notin V_t^+$. Since $A_t(r)-\Im
u(z)\in[0,\pi]$ by Lemma \ref{decreasing1}, the identity
\[\eta_\mu(z)=z\exp[-u(z)]=r\exp[iA_t(r)-u(z))]\] implies that
$\Im\eta_\mu(z)=0$ if and only if $A_t(r)=\Im u(z)$ or
$A_t(r)=\pi+\Im u(z)$. If $A_t(r)>0$ then $\Im u(z)=-A_t(r)/(t-1)<0$
by (\ref{Imu}). Further suppose that $A_t(r)=\pi+\Im u(z)$. This
gives $A_t(r)=(t-1)\pi/t$ and $\Im u(z)=-\pi/t$, and therefore
\[\arg\eta_\mu(z)=\arg z-\Im u(z)=\pi,\] which is a contradiction
since $\arg\eta_\mu(z)\in[A_t(r),\pi)$ by Proposition
\ref{chareta1}. This shows the necessity. Conversely, if $r\notin
V_t^+$, i.e., $A_t(r)=0$ by (\ref{Vt+1}) then $u(z)\in\mathbb{R}$,
and hence $\arg\eta_\mu(z)=\arg r-\Im u(r)=0$, and hence
$\Im\eta_\mu(z)=0$.

To verify the assertions (1) and (2), observe that
$\eta_{\mu^{\boxtimes
t}}(h_t(r))=(\eta_\mu\circ\omega_t)(\Phi_t(z))= \eta_\mu(z)$ and
\[\Im\left(\frac{1}{1-\eta_\mu(z)}\right)=\frac{\Im\eta_\mu(z)}{|1-\eta_\mu(z)|^2}\neq0\]
if and only if $r\in V_t^+$ by the above discussion. This implies
(1) since
\begin{align*}
\frac{d(\mu^{\boxtimes
t})^{\mathrm{ac}}}{dx}\left(\frac{1}{h_t(r)}\right)&=\frac{1}{\pi}\Im\left(\frac{h_t(r)}{1-\eta_{\mu^{\boxtimes
t}}(h_t(r))}\right) \\
&=\frac{h_t(r)}{\pi}\Im\left(\frac{1}{1-\eta_\mu(z)}\right) \\
&=\frac{h_t(r)\Im\eta_\mu(z)}{\pi|1-\eta_\mu(z)|^2}.
\end{align*}
by Stietljes inversion formula. The desired density for
$(\mu^{\boxtimes t})^{\mathrm{ac}}$ follows from the identities
\[A_t(r)-\Im
u(z)=A_t(r)+\frac{A_t(r)}{t-1}=\theta_t(r),\;\;\;\;\;r\in V_t^+,\]
\[\Re\eta_\mu(z)=l_t(r)\cos\theta_t(r)\;\;\;\;\;\mathrm{and}\;\;\;\;\;\Im\eta_\mu(z)=l_t(r)\sin\theta_t(r).\]
The analyticity of this density wherever it is positive follows from
Theorem \ref{Htprop}(4). Indeed, if $r\in V_t^+$ then
$\omega_t(h_t(r))=(\omega_t\circ\Phi_t)(re^{iA_t(r)})=re^{iA_t(r)}\in\mathbb{C}^+$,
which yields that $A_t$ is analytic in a neighborhood of $r$. To
prove the statement (3), it is enough to show that $V_t^+$ is a
decreasing set as $t$ increases. This will hold if we can show that
$g$ never has a local maximum in any open interval $I$ in $V_t^+$.
Indeed, the function $g$ is strictly convex on $I$ by Lemma
\ref{sconvex1}, whence (3) follows.
\end{pf} \qed

Recall that a point $x\in(0,\infty)$ is an atom for $\mu$ if and
only if $\eta_\mu(1/x)=1$ and $\eta_\mu'(1/x)$ is finite, in which
case we have
\[\eta_\mu'(1/x)=\frac{x}{\mu(\{x\})}.\]

\begin{prop} \label{prechar1}
If $r\in(0,\infty)$ and $t>1$ then the following
statements are equivalent.
\begin{enumerate} [$\qquad(1)$]
\item {$r\in(V_t)^c$ and $\eta_\mu(r)=1$;}
\item {$\eta_{\mu^{\boxtimes t}}(r^t)=1$;}
\item {$\mu(\{1/r\})\geq1-t^{-1}$.}
\end{enumerate}
If $(1)$-$(3)$ hold then
\begin{equation} \label{atomR}
1+\int_0^\infty\frac{r(s^2+1)}{(r-s)^2}\;d\rho(s)=\frac{1}{\mu(\{1/r\})}.
\end{equation}
\end{prop}

\begin{pf} The equivalence of (2) and (3) was proved in [\ref{BB2}].
We will show that (1) and (2) are equivalent. If (1) holds then
$A_t(r)=0$ by (\ref{Vt+1}), and hence $e^{u(r)}=r$ and
$h_t(r)=\Phi_t(r)=r\exp[(t-1)u(r)]=r^t$. This implies that
$\eta_{\mu^{\boxtimes t}}(r^t)=\eta_{\mu^{\boxtimes
t}}(h_t(r))=\eta_\mu(r)=1$ and (2) holds. Conversely, suppose that
(2) holds and $h_t(s)=r^t$ for some $s\in(0,\infty)$. Then we have
$\eta_\mu(se^{iA_t(s)})=\eta_{\mu^{\boxtimes t}}(r^t)=1$. By the
fact that $\arg\eta_\mu(z)\in[\arg z,\pi]$ for any
$z\in\mathbb{C}^+$, we must have $A_t(s)=0$ and $s\in(V_t^+)^c$.
Then the identities $se^{u(s)}=\eta_\mu(s)=1$ and
$s\exp[(t-1)u(s)]=\Phi_t(s)=r^t$ yield $r=s$, which implies (1).
Taking the Julia-Carath\'{e}dory derivative of
$\eta_\mu(z)=ze^{-u(z)}$ gives
$\eta_\mu'(z)=e^{-u(z)}-zu'(z)e^{-u(z)}$, and therefore
\[r\eta_\mu'(r)=re^{-u(r)}-ru'(r)re^{-u(r)}=1-ru'(r).\] On the other
hand, we have
\[u'(r)=\lim_{\epsilon\to0}\frac{u(re^{i\epsilon})-u(r)}{re^{i\epsilon}-r}=
\lim_{\epsilon\to0}\int_0^\infty\frac{-(s^2+1)d\rho(s)}{(re^{i\epsilon}-s)(r-s)}=
-\int_0^\infty\frac{(s^2+1)}{(r-s)^2}d\rho(s).\] Hence the equation
(\ref{atomR}) follows from the above discussions.
\end{pf} \qed

\begin{prop} \label{charatom1}
Let $r\in(0,\infty)$. Then the point $r^t$ is an atom of
$\mu^{\boxtimes t}$ if and only if $\mu(\{r\})>(t-1)/t$, in which
case,
\[\mu^{\boxtimes t}(\{r^t\})=t\mu(\{r\})-(t-1).\]
\end{prop}

\begin{pf} By Theorem \ref{Phit1prop}(4), it suffices to show the
necessity. Suppose that $r^t$ is an atom of $\mu^{\boxtimes t}$ and
let $s=1/r$. Then by the proof of Proposition \ref{prechar1} we see
that $h_t(s)=\Phi_t(s)=s^t$. Taking the Julia-Carath\'{e}dory
derivative of the identity $\eta_{\mu^{\boxtimes
t}}=\eta_\mu\circ\omega_t$ gives
\[\eta_{\mu^{\boxtimes
t}}(s^t)=\eta_\mu'(\omega_t(s^t))\omega_t'(s^t)=\frac{\eta_\mu'(s)}{\Phi'(s)},\]
where the fact that $\Phi_t'(\omega_t(s^t))\omega_t'(s^t)=1$. Since
$\eta_{\mu^{\boxtimes t}}(s^t)=r^t/\mu^{\boxtimes
t}(\{r^t\})<\infty$, we must have $\Phi_t'(s)>0$. Taking the
Julia-Carath\'{e}dory derivative of the identity
$\Phi_t(z)=z\exp[(t-1)u(z)]$ gives
\[\Phi_t'(z)=\exp[(t-1)u(z)]+(t-1)zu'(z)\exp[(t-1)u(z)],\] which
implies
\[s\Phi_t'(s)=\Phi_t(s)+(t-1)u'(s)\Phi_t(s)=\Phi_t(s)[1+(t-1)u'(s)].\]
This shows that $1+(t-1)su'(s)>0$ or, equivalently,
\[su'(s)<\frac{-1}{t-1}.\]
Since we have
\[su'(s)=1-\frac{1}{\mu(\{r\})}\] by the equation (\ref{atomR}), it
is easy to see that $\mu(\{r\})>1-t^{-1}$, as desired.
\end{pf} \qed

Combining Theorem \ref{density1} and Proposition \ref{charatom1}, we
have the following result.

\begin{thm} If $\mu\in\mathcal{M}_{\mathbb{R}_+}^\times$ then the following
statements hold.
\begin{enumerate} [$\qquad(1)$]
\item {The measure $\mu^{\boxtimes t}$, $t>1$, has at most countable many components in the support, which
consists of finitely many points $($atoms$)$ and countably many arcs
on $\mathbb{R}_+$.}
\item {The number of the components in $\mathrm{supp}(\mu^{\boxtimes t})$ is a decreasing function of $t$.}
\end{enumerate}
\end{thm}

\section{Support of the measure $\mu^{\boxtimes t}$ on $\mathbb{T}$}

Throughout this section, the measure
$\mu\in\mathcal{M}_\mathbb{T}^\times$ is fixed. By Proposition
\ref{chareta2}, there exists some analytic function
$u:\mathbb{D}\to\mathbb{C}$ with a nonnegative real part on
$\mathbb{D}$. A theorem of Herglotz yields that the function $u$ can
be expressed as
\begin{equation} \label{u1}
u(z)=i\alpha+\int_\mathbb{T}\frac{\zeta+z}{\zeta-z}\;d\rho(\zeta),\;\;\;\;\;z\in\mathbb{D}.
\end{equation} where $-\alpha=\arg|\eta_\mu'(0)|\in[-\pi,\pi]$ and
$\rho$ is a finite positive Borel measure on $\mathbb{T}$ satisfying
$\rho(\mathbb{T})=-\log|\eta_\mu'(0)|$. For $t>1$, let
\[\Phi_t(z)=z\exp[(t-1)u(z)],\;\;\;\;\;z\in\mathbb{D}.\]
By Proposition \ref{GIT2}, there exists a continuous function
$\omega_t:\overline{\mathbb{D}}\to\overline{\mathbb{D}}$ such that
$\omega_t|\mathbb{D}$ is analytic, $\omega_t(0)=0$, and
$\Phi_t(\omega_t(z))=z$ for all $z\in\mathbb{D}$. Let
\[\Omega_t=\{z\in\mathbb{D}:|\Phi_t(z)|<1\}.\] Then $\omega_t(\Phi_t(z))=z$,
$z\in\Omega_t$. The following proposition states another important
property of $\omega_t$.

\begin{prop} The function $\omega_t$ is the $\eta$-transform of some
$\boxtimes$-infinitely divisible measure in
$\mathcal{M}_\mathbb{T}^\times$.
\end{prop}

\begin{pf} The function $\omega_t$ satisfies the conditions in
Theorem \ref{chareta2}, and therefore it must be of the form
$\omega_t=\eta_{\nu_t}$ for some measure $\nu_t$ in
$\mathcal{M}_\mathbb{T}^\times$. Next, observe that the
$\Sigma$-transform of $\nu_t$ is given by
$\Sigma_{\nu_t}(z)=\Phi_t(z)/z=\exp[(t-1)u(z)]$. This implies the
desired result by [Theorem 7.5, \ref{HV1}].
\end{pf} \qed

Let $\mu^{\boxtimes t}$ be the unique measure in
$\mathcal{M}_\mathbb{T}$ satisfying
\begin{equation} \label{uniquemut}
\eta_{\mu^{\boxtimes
t}}(z)=\eta_\mu(\omega_t(z)),\;\;\;\;\;z\in\mathbb{D}.
\end{equation} Since
$\omega_t'(0)=1/\Phi_t'(\omega_t(0))=1/\Phi_t'(0)=\exp[-(t-1)u(0)]$,
we see that
\[\eta_{\mu^{\boxtimes
t}}'(0)=\eta_\mu'(0)\omega_t'(0)=\exp[-tu(0)],\] which particularly
shows that $\mu^{\boxtimes t}\in\mathcal{M}_\mathbb{T}^\times$.

In the rest of this section, we will investigate the support of
$\mu^{\boxtimes t}$, which is the measure satisfying the requirement
(\ref{uniquemut}). Our analysis will be based on the following
functions $g:[-\pi,\pi]\to\mathbb{R}^+\cup\{+\infty\}$ and
$R_t:[-\pi,\pi]\to[0,1]$, which are defined as
\begin{equation} \label{g2}
g(\theta)=\int_{-\pi}^\pi\frac{d\rho(e^{i\phi})}{1-\cos(\theta-\phi)}
\end{equation} and
\begin{equation} \label{Rt}
R_t(\theta)=\sup\left\{r\in(0,1):\int_{-\pi}^\pi
T(r,\theta-\phi)\;d\rho(e^{i\phi})<\frac{1}{t-1}\right\},
\end{equation} respectively, where
\begin{equation} \label{Ttheta}
T(r,\theta)=\frac{r^2-1}{\log r}\frac{1}{1-2r\cos\theta+r^2}
\end{equation} is a continuous
function from $(0,1)\times[-\pi,\pi]$ to $\mathbb{R}^+$. The following set, associated with the
function $g$,
\[V_t^+=\left\{\theta\in[-\pi,\pi]:g(\theta)>\frac{1}{t-1}\right\},\]
will also play an important role in the study of the support of
$\mu^{\boxtimes t}$.

The following lemmas provide some basic properties about the
functions and set introduced above.

\begin{lem} \label{decreasing2}
For any $\theta\in[-\pi,\pi]$, the function $T(\cdot,\theta)$
defined as in $(\ref{Ttheta})$ is strictly increasing on $(0,1)$.
Consequently,
\[T(1,\theta)=\lim_{r\uparrow1}T(r,\theta)=\frac{1}{1-\cos\theta}\]
for $\theta\neq0$ and $T(1,0)=\lim_{r\uparrow1}T(r,0)=+\infty$.
\end{lem}

\begin{pf} First note that the function
\[r\mapsto\frac{(r+1)^2}{1-2r\cos\theta+r^2},\;\;\;\;\;r\in(0,1),\] is increasing since
it has a non-negative derivative on $(0,1)$. Next, by some simple
manipulations we see that the function
\[f(x)=\frac{e^x-1}{x(e^x+1)}\] has a strictly negative
derivative on $(0,\infty)$. This implies that the function
\[r\mapsto f(-\log r)=\frac{r-1}{(r+1)\log r}\] is strictly increasing on $(0,1)$.
Then the desired conclusion follows from the expression
\[T(r,\theta)=\frac{r-1}{(r+1)\log
r}\frac{(r+1)^2}{1-2r\cos\theta+r^2}.\] This completes the proof.
\end{pf} \qed

\begin{lem} Let $t>1$ and $T(r,\theta)$ be the function defined in $(\ref{Ttheta})$ and
$\mathrm{Lemma}$ $\ref{decreasing2}$ on $(0,1]\times[-\pi,\pi]$.
Then
\begin{equation} \label{Omegat2}
\Omega_t=\{re^{i\theta}:r\in[0,R_t(\theta))\;\;\mathrm{and}\;\;\theta\in[-\pi,\pi]\}
\end{equation} and
\begin{equation} \label{partialomegat2}
\partial\Omega_t=\{R_t(\theta)e^{i\theta}:\theta\in[-\pi,\pi]\}.
\end{equation} For any $\theta\in[-\pi,\pi]$ we have
\begin{equation} \label{Rtran}
R_t(\theta)\in(0,1]
\end{equation}
and
\begin{equation} \label{Vt+2}
V_t^+=\{\theta\in[-\pi,\pi]:R_t(\theta)<1\}.
\end{equation} Moreover, for any
$\theta\in[-\pi,\pi]$, we have
\begin{equation} \label{integralRt}
\int_{-\pi}^\pi
T(R_t(\theta),\theta+\phi)\;d\rho(e^{i\phi})\leq\frac{1}{t-1},
\end{equation}
where the equality holds if $\theta\in V_t^+$.
\end{lem}

\begin{pf} First observe that a point $re^{i\theta}$ belongs to
$\Omega_t\backslash\{0\}$ if and only if
\[g(r,\theta):=\int_{-\pi}^\pi T(r,\theta-\phi)\;d\rho(e^{i\phi})<1/(t-1).\] Indeed, simple computations give
\[\log|\Phi_t(re^{i\theta})|=(\log r)[1-(t-1)g(r,\theta)].\] This
implies that $re^{i\theta}\in\mathbb{D}\backslash\{0\}$ belongs to
$\Omega_t$ if and only if $\log|\Phi_t(re^{i\theta})|<0$, which
happens if and only if $(\log r)[1-(t-1)g(r,\theta)]<0$ or,
equivalently, $g(r,\theta)<1/(t-1)$. By Lemma \ref{decreasing2} and
the definition of the function $R_t$ we see that the identity
(\ref{Omegat1}) holds. The identity (\ref{partialomegat2}) follows
directly from (\ref{Omegat2}).

For any $\theta\in[-\pi,\pi]$ and small $\epsilon>0$, by the fact
(\ref{Omegat2}) we have $g(R_t(\theta)-\epsilon,\theta)<1/(t-1)$,
which gives (\ref{integralRt}) by letting $\epsilon\to0$ and
monotone convergence theorem. This also shows that $\theta\in V_t^+$
if $R_t(\theta)<1$. If $R_t(\theta)=1$ then
$g(r,e^{i\theta})\leq1/(t-1)$ for $r\in(0,1)$, whence
$g(\theta)=\lim_{r\uparrow1^-}g(r,e^{i\theta})\leq1/(t-1)$. Hence
the identity (\ref{Vt+2}) holds. By the above discussion, we see
that the equality in (\ref{integralRt}) holds if $\theta\in V_t^+$.
\end{pf} \qed

Due to the preceding lemma, we have the the following proposition
which generalizes the statement (3) in Theorem \ref{GIT2}.

\begin{prop} For any $z\in\Omega_t$, the line segment joining the origin and
$z$ is contained in $\Omega_t$. Consequently, the set $\Omega_t$
consists of one connected component.
\end{prop}

\begin{prop} \label{Lip2}
The function $u(z)$ has a continuous extension to
$\overline{\Omega}_t$. Moreover, the extension is Lipschitz
continuous on $\overline{\Omega}_t$ and can be expressed as
\begin{equation} \label{extension2}
u(z)=i\alpha+\int_\mathbb{T}\frac{\zeta+z}{\zeta-z}\;d\rho(\zeta),\;\;\;\;\;\zeta\in\overline{\Omega_t}.
\end{equation}
\end{prop}

\begin{pf} We first show that the integral in (\ref{extension2})
converges for $z\in\overline{\Omega}_t$. It suffices to consider the
case $z=e^{i\theta}\in\partial\Omega_t\cap\mathbb{T}$. Note that
$z_r:=rz\in\Omega_t$, $0<r<1$, and
\[\int_\mathbb{T}\frac{\zeta+z_r}{\zeta-z_r}\;d\rho(\zeta)=
\int_\mathbb{T}\frac{1-r^2}{|\zeta-z_r|^2}\;d\rho(\zeta)+
2i\int_{-\pi}^\pi\frac{r\sin(\theta-\phi)d\rho(e^{i\phi})}{1-2r\cos(\theta-\phi)+r^2}.\]
The first integral on the right hand of the above expression tends
to zero as $r\to1^-$ since
\[\int_\mathbb{T}\frac{1-r^2}{|\zeta-z_r|^2}\;d\rho(\zeta)\leq\frac{-\log r}{t-1}\]
by (\ref{integralRt}). By the fact that the mapping
\[r\mapsto\frac{r}{1-2r\cos\phi+r^2}\] is strictly increasing on
$(0,1)$ for any $\phi$, we conclude that the integral
\[\int_{-\pi}^\pi\frac{r\sin(\theta-\phi)}{1-2r\cos(\theta-\phi)+r^2}\;d\rho(e^{i\phi})\to
\frac{1}{2}\int_{-\pi}^\pi\frac{\sin(\theta-\phi)}{1-\cos(\theta-\phi)}\;d\rho(e^{i\phi})\]
converges as $r\to1^-$ by monotone convergence theorem. Since
$g(\theta)\leq1/(t-1)$ by (\ref{integralRt}), it follows that $u$
has a continuous extension to $\overline{\Omega}_t$. Next, we show
that this extension is Lipschitz continuous. Notice that for any
different points $z_1$ and $z_2$ in $\Omega_t$ we have
\begin{align*}
\frac{|u(z_1)-u(z_2)|}{|z_1-z_2|}&=2\left|\int_\mathbb{T}\frac{\zeta
d\rho(\zeta)}{(\zeta-z_1)(\zeta-z_2)}\right| \\
&\leq2\left(\int_\mathbb{T}\frac{d\rho(\zeta)}{|\zeta-z_1|^2}\right)^{1/2}
\left(\int_\mathbb{T}\frac{d\rho(\zeta)}{|\zeta-z_2|^2}\right)^{1/2},
\end{align*} where the H\"{o}lder inequality is used in the the last inequality.
Let $\epsilon$ be small enough so that the set
$\overline{\Omega}_{t,\epsilon}:=\overline{\Omega}_t\cap\{z:|z|\geq\epsilon\}\neq\emptyset$
(note that $V_t^+$ is an open set containing the origin). Then the
fact that the function $r\mapsto\log r/(r^2-1)$ is decreasing on
$(0,1)$ shows that
\[\int_\mathbb{T}\frac{d\rho(\zeta)}{|\zeta-z|^2}\leq
\frac{\log|z|}{|z|^2-1}\frac{1}{t-1}\leq\frac{\log\epsilon}{\epsilon^2-1}\frac{1}{t-1},\;\;\;\;\;z\in
\overline{\Omega}_{t,\epsilon}.\] On the other hand, if
$z\in\Omega_t$ with $|z|\leq\epsilon$ then
\[\int_\mathbb{T}\frac{d\rho(\zeta)}{|\zeta-z|^2}\leq
\int_\mathbb{T}\frac{d\rho(\zeta)}{(1-\epsilon)^2}=(1-\epsilon)^{-2}\rho(\mathbb{T}).\]
The above discussions yield the desired result.
\end{pf} \qed

For any $-\pi\leq a<b\leq\pi$, let
$A_{a,b}=\{e^{i\theta}:a<\theta<b\}$ be an arc contained in
$\mathbb{T}$.

\begin{lem} \label{sconvex2}
If the function $g$ is bounded on some open interval $(a,b)$ then
$\rho(A_{a,b})=0$ and $g_t$ is strictly convex on $(a,b)$. In
particular, this is true if $(a,b)$ is contained in $(V_t^+)^c$.
\end{lem}

\begin{pf} Suppose that $g$ is bounded by $M$ on $(a,b)$. Let
$(c,d)$ be any subinterval of $(a,b)$. If $\theta=(c+d)/2$ then we
have
\[M\geq\int_{-\pi}^\pi\frac{d\rho(e^{i\phi})}{1-\cos(\theta-\phi)}\geq
\int_c^d\frac{d\rho(e^{i\phi})}{1-\cos(\theta-\phi)}=
\int_c^d\frac{d\rho(e^{i\phi})}{2\sin^2\left(\frac{\theta-\phi}{2}\right)}.\]
Since $\sin^2\phi\leq\phi^2$ for any $\phi\in[-\pi,\pi]$, we deduce
that
\[M\geq2\int_a^b\frac{d\rho(e^{i\phi})}{(\theta-\phi)^2}\geq\frac{8\rho(A_{c,d})}{(d-c)^2}\]
or, equivalently,
\[\frac{\rho(A_{c,d})}{d-c}\leq8M(d-c).\] Since the above inequality holds for any
subinterval contained in $(a,b)$, we conclude that the desired
result holds. The strict positivity of the second order derivative
\[g''(\theta)=\frac{3}{4}\int_{[-\pi,\pi]\backslash(a,b)}\frac{d\rho(e^{i\theta})}
{\sin^4\left(\frac{\theta-\phi}{2}\right)}\] on $I$ yield the second
assertion.
\end{pf} \qed

Observe that the function $\theta\mapsto R_t(\theta)e^{i\theta}$ is
a homeomorphism from $[-\pi,\pi]$ onto $\partial\Omega_t$. Since
$\Phi_t$ has a continuous extension to $\overline{\Omega}_t$, we
deduce that the function
\[h_t(e^{i\theta})=\Phi_t(R_t(\theta)e^{i\theta})\] is a homeomorphism
of $\mathbb{T}$. Now, we are in a position to state the main theorem
of this section.

\begin{thm} \label{density2}
Suppose that $\mu$ is a measure in $\mathcal{M}_\mathbb{T}^\times$
and $t>1$. Let
\[S_t=\{\overline{h_t(e^{i\theta})}:\theta\in V_t^+\}.\] Then the following statements hold.
\begin{enumerate} [$\qquad(1)$]
\item {The measure $(\mu^{\boxtimes t})^{\mathrm{ac}}$ is concentrated on the closure of $S_t$.}
\item {The density of $(\mu^{\boxtimes t})^{\mathrm{ac}}$ is analytic on the set
$S_t$ and is given by
\[\frac{d\mu^{\boxtimes t}}{d\zeta}\left(\overline{h_t(e^{i\theta})}\right)=
\frac{1}{2\pi}\frac{1-l_t^2(\theta)}{1-2l_t(\theta)\cos\alpha(\theta)+l_t^2(\theta)},\;\;\;\;\;\theta\in
V_t^+,\] where
\[l_t(\theta)=R_t^{\frac{t}{t-1}}(\theta)\]
\[\alpha(\theta)=\theta-\Im u(R_t(\theta)e^{i\theta}).\]}
\item {The number of components in $\mathrm{supp}(\mu^{\boxtimes t})^{\mathrm{ac}}$ is a decreasing function of $t$.}
\end{enumerate}
\end{thm}

\begin{pf} Let $z=R_t(\theta)e^{i\theta}$. We claim that
$|\eta_\mu(z)|=1$ if and only if $\theta\not\in V_t^+$. Observe that
if $\theta\in(V_t^+)^c$ then $R_t(\theta)=1$ and $\Re
u(e^{i\theta})=0$, which gives
\[|\eta_\mu(e^{i\theta})|=|e^{i\theta}e^{-\Re u(e^{i\theta})}|=1.\]
If $\theta\in V_t^+$ then
\[\Re u(z)=\frac{-\log R_t(\theta)}{t-1}\]
by (\ref{integralRt}), whence we have
\[|\eta_\mu(z)|=R_t^{\frac{t}{t-1}}(\theta)\] and
\[\Re\eta_\mu(z)=R_t^{\frac{t}{t-1}}(\theta)\cos(\theta-\Im u(z)).\]
Since $\eta_{\mu^{\boxtimes t}}$ is continuous on $\mathbb{T}$, it
follows that $(\mu^{\boxtimes t})^\mathrm{ac}$ is concentrated on
the closure of the set of points $\zeta\in\mathbb{T}$ such that
\[\Re\left(\frac{1+\eta_{\mu^{\boxtimes t}}(\zeta)}{1-\eta_{\mu^{\boxtimes t}}(\zeta)}\right)
\] is finite and nonzero. Now, note that we have $\eta_{\mu^{\boxtimes
t}}(h_t(e^{i\theta}))=\eta_\mu(z)$ and
\[\Re\left(\frac{1+\eta_{\mu^{\boxtimes t}}(h_t(e^{i\theta}))}{1-\eta_{\mu^{\boxtimes t}}(h_t(e^{i\theta}))}\right)=
\Re\left(\frac{1+\eta_\mu(z)}{1-\eta_\mu(z)}\right)=\frac{1-|\eta_\mu(z)|^2}{1-2\Re\eta_\mu(z)+|\eta_\mu(z)|^2}.\]
Since
\[\frac{d(\mu^{\boxtimes t})^{\mathrm{ac}}}{d\zeta}\left(\overline{h_t(e^{i\theta})}\right)
=\Re\left(\frac{1+\eta_{\mu^{\boxtimes
t}}(h_t(e^{i\theta}))}{1-\eta_{\mu^{\boxtimes
t}}(h_t(e^{i\theta}))}\right),\] the desired results in (1) and (2)
follow. To verify the statement (3), it suffices to show that the
number of components in $V_t^+$ is nonincreasing as $t$ increases.
This will hold if we show that $g$ never has a local maximum in any
open interval $(a,b)$ in $V_t^+$. Indeed, the function $g$ is
strictly convex on such an interval by Lemma \ref{sconvex2}, whence
(3) follows.
\end{pf} \qed

For the rest of this section, we discuss the atoms of
$\mu^{\boxtimes t}$.

\begin{prop} \label{prechar2}
If $\theta\in[-\pi,\pi]$ and $t>1$ then the following statements are
equivalent.
\begin{enumerate} [$\qquad(1)$]
\item {$\theta\in(V_t^+)^c$ and $\eta_\mu(e^{i\theta})=1$;}
\item {$\eta_{\mu^{\boxtimes t}}(h_t(e^{i\theta}))=1$;}
\item {$\mu(\{e^{-i\theta}\})\geq1-t^{-1}$.}
\end{enumerate}
If $(1)$-$(3)$ hold then
\begin{equation} \label{atomT}
1+\int_{-\pi}^\pi\frac{d\rho(e^{i\phi})}{1-\cos(\theta-\phi)}=\frac{1}{\mu(\{e^{i\theta}\})}.
\end{equation}
In addition, if $t\geq2$ then the condition $\eta_{\mu^{\boxtimes
t}}(e^{it\theta})=1$ with $|t\theta|<\pi$ is equivalent to
conditions $(1)$-$(3)$.
\end{prop}

\begin{pf} The equivalence of (2) and (3) was proved in [\ref{BB2}].
We will show that (1) and (2) are equivalent. If (1) holds then
$R_t(\theta)=1$ by (\ref{Rtran}), and hence $\eta_{\mu^{\boxtimes
t}}(h_t(e^{i\theta}))=\eta_\mu(\omega_t(\Phi_t(e^{i\theta})))=\eta_\mu(e^{i\theta})=1$,
which yields (2). Conversely, the condition
$\eta_\mu(R_t(\theta)e^{i\theta})=\eta_{\mu^{\boxtimes
t}}(h_t(e^{i\theta}))=1$ in (2) along with the fact that
$|\eta_\mu(z)|\leq|z|$ for any $z\in\mathbb{D}$ shows
$R_t(\theta)=1$, and therefore (1) holds. This also particularly
implies that
\begin{equation} \label{angle}
\exp[u(e^{i\theta})]=e^{i\theta}.
\end{equation}

Next, suppose that $\eta_\mu(e^{i\theta})=1$ and
$e^{i\theta}\eta_\mu'(e^{i\theta})<\infty$, where
$\eta_\mu'(e^{i\theta})$ is the Julia-Carath\'{e}dory derivative of
$\eta_\mu$ at $e^{i\theta}$. Taking the Julia-Carath\'{e}dory
derivative of the identity $\eta_\mu(z)=z\exp[-u(z)]$ gives
\[\eta_\mu'(z)=\exp[-u(z)]-zu'(z)\exp[-u(z)],\] and
therefore we have
\begin{align*}
e^{i\theta}\eta_\mu'(e^{i\theta})&=e^{i\theta}\exp[-u(e^{i\theta})]-e^{i\theta}u'(e^{i\theta})e^{i\theta}
\exp[-u(e^{i\theta})] \\
&=1-e^{i\theta}u'(e^{i\theta}).
\end{align*} On the other hand, using
the identity
\[\frac{u(re^{i\theta})-u(e^{i\theta})}{(r-1)e^{i\theta}}=\int_\mathbb{T}\frac{2\xi d\rho(\xi)}{(\xi-re^{i\theta})
(\xi-e^{i\theta})},\;\;\;\;\;0<r<1,\] it is easy to see that the
Julia-Carath\'{e}dory derivative $u'(e^{i\theta})$ is given by
\[u'(e^{i\theta})=\int_\mathbb{T}\frac{2\xi d\rho(\xi)}{(\xi-e^{i\theta})^2}.\]
Then
\begin{align*}
e^{i\theta}u'(e^{i\theta})&=\Re\left(\int_\mathbb{T}\frac{2e^{i\theta}\xi}{(\xi-e^{i\theta})^2}\;d\rho(\xi)\right) \\
&=\int_\mathbb{T}\Re\left(\frac{2e^{i\theta}\xi}{(\xi-e^{i\theta})^2}\right)d\rho(\xi) \\
&=-\int_{-\pi}^\pi\frac{d\rho(e^{i\phi})}{1-\cos(\theta-\phi)},
\end{align*} which gives the identity (\ref{atomT}).

Finally, if (2) holds and $t\geq2$ then
\[|u(e^{i\theta})|=\left|\int_\mathbb{T}\frac{2\Im(e^{i\theta}\overline{\xi})}
{|\xi-e^{i\theta}|^2}\;d\rho(\xi)\right|\leq\int_{-\pi}^{\pi}\frac{d\rho(e^{i\phi})}
{1-\cos(\theta-\phi)}\leq\frac{1}{t-1}<\pi,\] which shows that
$u(e^{i\theta})=i\theta$ by the equation (\ref{angle}). Therefore
\[h_t(e^{i\theta})=\Phi_t(e^{i\theta})=e^{i\theta}\exp[(t-1)u(e^{i\theta})]=e^{it\theta}\]
and $|t\theta|=|tu(e^{i\theta})|\leq t/(t-1)<\pi$. Now, suppose
$\eta_{\mu^{\boxtimes t}}(e^{it\theta})=1$ with $|t\theta|<\pi$ and
$h_t(e^{i\phi})=e^{it\theta}$ for some $\phi\in[-\pi,\pi]$. Then the
preceding argument indicates that $u(e^{i\phi})=i\phi$,
$|\phi|\leq1/(t-1)$ and $h_t(e^{i\phi})=e^{it\phi}$. Since
$e^{it\theta}=e^{it\phi}$ and $|t\theta|,|t\phi|\leq\pi$, we must
have $\theta=\phi$. This completes the proof.
\end{pf} \qed

\begin{prop} \label{charatom2}
A point $1/\zeta$ is an atom of $\mu^{\boxtimes t}$ if and only if
$\mu(\{1/\omega_t(\zeta)\})>1-t^{-1}$, in which case we have
\[\mu^{\boxtimes t}(\{1/\zeta\})=t\mu(\{1/\omega_t(\zeta)\})-(t-1).\]
\end{prop}

\begin{pf} By Theorem \ref{Phit2prop}(5), it suffices to show the necessity. Suppose that
$1/\zeta$ is an atom of $\mu^{\boxtimes t}$ and let
$h_t(e^{i\theta})=\xi$ for some $\theta\in[-\pi,\pi]$. Then by
Proposition \ref{prechar1} we see that $\Phi_t(e^{i\theta})=\xi$.
Taking the Julia-Carath\'{e}dory derivative of the identity
$\eta_{\mu^{\boxtimes t}}=\eta_\mu\circ\omega_t$ gives
\[\xi\eta_{\mu^{\boxtimes t}}'(\xi)=\xi\eta_\mu'(\omega_t(\xi))\omega_t'(\xi)=
e^{i\theta}\eta_\mu'(e^{i\theta})\frac{\xi}{e^{i\theta}\Phi_t'
(e^{i\theta})},\] where the fact that
$\Phi_t'(\omega_t(\xi))\omega_t'(\xi)=1$. Since
$\xi\eta_{\mu^{\boxtimes t}}'(\xi)=1/\mu^{\boxtimes
t}(\{1/\xi\})<\infty$, we must have
\begin{equation} \label{finite2}
\frac{\xi}{e^{i\theta}\Phi_t' (e^{i\theta})}<\infty.
\end{equation}
Taking the Julia-Carath\'{e}dory derivative of the identity
$\Phi_t(z)=z\exp[(t-1)u(z)]$ gives
\[\Phi_t'(z)=\exp[(t-1)u(z)]+(t-1)u'(z)z\exp[(t-1)u(z)],\] which
implies
\[e^{i\theta}\Phi_t'(e^{i\theta})=\Phi_t(e^{i\theta})+(t-1)e^{i\theta}u'(e^{i\theta})\Phi_t(e^{i\theta})=
\xi[1+(t-1)e^{i\theta}u'(e^{i\theta})].\] Using the above identity,
the condition (\ref{finite2}) is equivalent to
\[e^{i\theta}u'(e^{i\theta})>\frac{-1}{t-1}.\] Since we have
\[e^{i\theta}u'(e^{i\theta})=1-\frac{1}{\mu(\{e^{-i\theta}\})}\] by
Proposition \ref{prechar1}, it is easy to see that
$\mu(\{e^{-i\theta}\})>1-t^{-1}$, as desired.
\end{pf} \qed

Combining Proposition \ref{charatom2} and Theorem \ref{density2}, we
have the following result.

\begin{thm} If $\mu\in\mathcal{M}_\mathbb{T}^\times$ then the following
statements hold.
\begin{enumerate} [$\qquad(1)$]
\item {The measure $\mu^{\boxtimes t}$, $t>1$, has at most countable many components in the support, which
consists of finitely many points $($atoms$)$ and countably many arcs
on $\mathbb{T}$.}
\item {The number of the components in $\mathrm{supp}(\mu^{\boxtimes t})$ is a decreasing function of $t$.}
\end{enumerate}
\end{thm}

\section*{Acknowledgments} The authors wish to thank their advisor,
Professor Hari Bercovici, for his generosity, encouragement, and
invaluable discussion during the course of the investigation.


\begin{thebibliography}{99}

\bibitem{} \label{moment} N. I. Achieser, {\it The classical moment problem}, in
Russian, Fizmatgiz, Moscow, 1961.

\bibitem{} \label{Japan} O. Arizmendi, T. Hasebe, Semigroups related to additive and multiplicative, free and Boolean
convolutions. Arxiv:1105.3344v3.

\bibitem{} \label{BB1} S.T. Belinschi, H. Bercovici, Atoms and regularity for measures in a
partially defined free convolution semigroup, {\it Math. Z.}
{\bf248} (4) 665-674 (2004).

\bibitem{} \label{BB2} S.T. Belinschi, H. Bercovici, Partially defined semigroups relative
to multiplicative free convolution, {\it Int. Math. Res. Not.}
{\bf2} 65-101 (2005).

\bibitem{} \label{BB3} S.T. Belinschi, H. Bercovici, A new approach to subordination results
in free probability, {\it J. D'analyse Math.} {\bf101} 357-365
(2007).

\bibitem{} \label{BN1} S.T. Belinschi, A. Nica, On a remarkable semigroup of
homomorphisms with respect to free multiplicative convolution, {\it
Indiana. Univ. Math J.} {\bf57} (4) 1679-1713 (2008).

\bibitem{} \label{BN2} S.T. Belinschi, A. Nica, Free Brownian motion and evolution towards $\boxplus$-infinitely
divisibility for $k$-tuples, {\it Int. J. Math.} {\bf20} (3) 309-338
(2009).

\bibitem{} \label{BP} H. Bercovici, V. Pata, Stable laws and domains of attraction in
free probability theory, {\it Ann. of Math.}, {\bf149}, 1023-1060
(1999).

\bibitem{} \label{HV1} H. Bercovici, D. Voiculescu, L\'{e}vy-Hin\v{c}in type theorems for
multiplicatrive and additive free convolution {\it Pacific Journal
of Mathematics} {\bf153} No.2 217-248 (1992).

\bibitem{} \label{HV2} H. Bercovici, D. Voiculescu, Free Convolutions of measures with
unbounded support, {\it Indiana Univ. Math. J.} {\bf42} (3) 733-773
(1993).

\bibitem{} \label{HV3} H. Bercovici, D. Voiculescu, Superconvergence to
the central limit and failure of the Cram\'{e}r theorem for free
random variables, {\it Probab. Theory Relat. Fields} {\bf103}
215-222 (1995).

\bibitem{} \label{BV4} H. Bercovici, D. Voiculescu, Regularity questions for free convolution, in:
Nonselfadjoint Operator Algebras, Operator Theory, and Related
Topics, in: Oper. Theory Adv. Appl., vol. 104, Birkhauser, Basel,
1998, pp. 37-47.

\bibitem{} \label{Biane1} P. Biane, On the free convolution with a semi-circular distribution,
{\it Indiana Univ. Math. J.} {\bf46} (3) 705-718 (1997).

\bibitem{} \label{Biane2} P. Biane, Segal-Bargmann transform, functional calculus on
matrix spaces and the theory of semi-circular and circular systems,
{\it J. Funct. Anal.} {\bf144} no. 1, 232-286 (1997).

\bibitem{} \label{Biane3} P. Biane, Processes with free increments, {\it Math. Z.} {\bf227} (1) 143-174 (1998).

\bibitem{} \label{G1} G. P. Chistyakov, F. G\"{o}tze, Limit theorems in free
probability theory. I, {\it Ann. Probab.} {\bf36} No.1 54-90 (2008).

\bibitem{} \label{G4} G. P. Chistyakov, F. G\"{o}tze, The arithmetic of distributions
in free probability theory, {\it Cent. Eur. J. Math.} {\bf9} No.5
997-1050 (2011).

\bibitem{} \label{G2} G. P. Chistyakov, F. G\"{o}tze, Asymptotic Expansions in the CLT in Free
Probability. ArXiv: 1109.4844.

\bibitem{} \label{G3} G. P. Chistyakov, F. G\"{o}tze, Rate of Convergence in the entropic free CLT. ArXiv: 1112.5087.

\bibitem{} \label{Maa} H. Maassen, Addition of freely independent random variables, {\it
J. Funct. Anal.} {\bf106} 409-438 (1992).

\bibitem{} \label{Nica} A. Nica, Multi-variable subordination distributions for free additive convolution,
{\it J. Funct. Anal.} {\bf257} 428-463 (2009).

\bibitem{} \label{NS} A. Nica, R. Speicher, On the multiplication of free $N$-tuples of noncommutative
random variables. {\it Amer. J. Math.} {\bf118}(4), 799-837 (1996).

\bibitem{} \label{Boolean} R. Speicher, R. Woroudi, Boolean convolution, in
free probability theory, Ed. D. Voiculescu, {\it Fields. Inst.
Commun.} {\bf12} 267-280 (1997).

\bibitem{} \label{V1} D.V. Voiculescu, Addition of certain non-commuting random
variables, {\it J. Funct. Anal.} {\bf66} 323-346(1986).

\bibitem{} \label{V2} D.V. Voiculescu, The analogues of entropy and of Fisher's
information measure in free probability theory I, {\it Comm. Math.
Phys.} {\bf 155} (1) 411-440 (1993).

\bibitem{} \label{V3} D.V. Voiculescu, The coalgebra of the free difference
quotient and free probability, {\it Internat. Math. Res. Notices}
{\bf2} 79-106 (2000).

\bibitem{} \label{V4} D.V. Voiculescu, K.J. Dykema, A. Nica, Free Random Variables. CRM Monograph
Series, Vol. 1 Am. Math. Soc. Providence, RI, (1992).

\bibitem{} \label{P1} P. Zhong, Free Bronian motion and free
convolution semigroups: multiplicative case, arXiv:1210.6090.

\bibitem{} \label{P2} P. Zhong, On the free convolution wth a free multiplicative anaogue of the normal
distribution, arXiv:1211.3160.

\end{thebibliography}
\end{document}